\documentclass[10pt,oneside,a4paper]{amsart}
\usepackage{latexsym}
\usepackage{amsfonts,amsmath,amssymb,indentfirst}
\newcommand{\beq}{\begin{equation}}
\newcommand{\eeq}{\end{equation}}
\newcommand{\bdism}{\begin{displaymath}}
\newcommand{\edism}{\end{displaymath}}

\newcommand{\setE}{\mathbb E}
\newcommand{\setH}{\mathbb H}

\newcommand{\setZ}{\mathbb Z}

\newcommand{\setR}{\mathbb R}

\newcommand{\setS}{\mathbb S}
\newtheorem{theorem}{Theorem}[section]
\newtheorem{proposition}[theorem]{Proposition}

\newtheorem{remark}[theorem]{Remark}

\author{\scshape Luca Fabrizio Di Cerbo}
\title{\bf A gap property for the growth of closed 3-manifold groups}
\begin{document}

\thanks{*Supported in part by a Marie Curie  Fellowship and a Renaissance Technologies Fellowship}
\address{Department of Mathematics, SUNY, Stony Brook, NY 11794-3651,
USA} \email{luca@math.sunysb.edu}
\begin{abstract}
We provide a lower bound for the uniform exponential growth rate of closed nonflat nonpositively curved 3-manifold groups. A detailed study of the uniform exponential growth rate of closed 3-manifold groups is also presented.
\end{abstract}
\maketitle

\section{Introduction}
\pagenumbering{arabic}

Let $G$ be a finitely generated group and let $S$ be a finite generating set. The growth function of the pair $(G, S)$ associates to a integer $k$ the number $\gamma_{S}(k)$ of elements $g\in G$ with word length $l_{S}(k)\leq k$. The exponential growth rate of the pair $(G, S)$ is defined as
\begin{align}
\omega(G, S):= \lim_{k\rightarrow \infty}\sqrt[k]{\gamma_{S}(k)}
\end{align}
and the group $G$ is said to be of exponential growth if $\omega(G, S)> 1$. The uniform exponential growth rate of the pair $(G, S)$ is defined as
\begin{align}
\omega(G):=\inf_{S}\omega(G, S)
\end{align}
where $S$ runs through the set of finite generating subsets of $G$, and the group $G$ is said to be of uniform exponential growth if $\omega(G)>1$. We refer to \cite{Harpe} for the necessary background on the growth of groups.

The study of exponential growth rate and uniform exponential growth rate of the fundamental group of a Riemannian manifold $(M, g)$ is of interest in dynamical system theory and Riemannian geometry. As an example we recall the following well known result of J. Milnor, see \cite{Milnor1}.
\begin{theorem}[Milnor]\label{Milnor}
If $(M, g)$ is a closed manifold with strictly negative sectional curvature, then $\pi_{1}(M)$ is of exponential growth.
\end{theorem}
Here is an interesting generalization due to A. Avez, see \cite{Avez}.
\begin{theorem}[Avez]\label{Avez}
If $(M, g)$ is a closed manifold with nonpositive sectional curvature, then either $\pi_{1}(M)$ is of exponential growth or $M$ is flat.
\end{theorem}
In dimension two we have the following "uniform" strengthening of theorem \ref{Avez}.
\begin{proposition}\label{Harpe1}
If $(M^{2}, g)$ is a closed orientable surface with nonpositive sectional curvature, then either $M^{2}$ is flat or $\omega(\pi_{1}(M^{2}))\geq 5$.
\end{proposition}
This proposition is a consequence of the classical uniformization theorem for 2-dimensional closed Riemannian manifolds and of the estimate $\omega(\pi_{1}(F_{g}))\geq 4g-3$ for which we refer to \cite{Harpe}, where with $F_{g}$ we indicate the surface of genus $g$. Note that the weaker estimate $\omega(\pi_{1}(F_{g}))\geq 2g-1$ is readily derived.

In section section \ref{12} we address the problem of a suitable generalization of proposition \ref{Harpe1} to the class of closed nonpositively curved three manifolds.

In section \ref{12} we discuss naturally related problems in the totality of closed 3-manifold groups.

\section{On the growth of closed nonpositively curved three manifold groups}\label{12}
In this section we provide a proof of the following theorem.
\begin{theorem}\label{Di Cerbo}
 There exists a universal constant $C>1$ with the following property:\\
 if $(M^{3}, g)$ is any nonpositively  curved closed manifold then it is either flat or $\omega(\pi_{1}(M^{3}))\geq C$.
\end{theorem}
The complete proof of theorem \ref{Di Cerbo} relies on the recently proved uniformization theorem for closed 3-manifolds for which we refer to the original papers by Perelman or to \cite{Morgan}, \cite{france} and the bibliography there.

Let $M^{3}$ be a compact manifold of nonpositive sectional curvature. A theorem of Hadamard implies that $M^{3}$ must be covered by $\setR^{3}$. The Alexander theorem implies that $M^{3}$ must be irreducible, see \cite{Hempel}.
Once we reduced our study to the class of irreducible closed 3-manifolds we can invoke Perelman's uniformization result.
\begin{theorem}[Perelman]\label{Perelman}
Let $M^{3}$ be a closed, irreducible manifold. Then $M^{3}$ is hyperbolic, Seifert fibred, or its JSJ decomposition is nonempty.
\end{theorem}
This powerful classification result divides the proof of theorem \ref{Di Cerbo} into three distinct parts.
\subsection*{Hyperbolic 3-manifolds}
A uniform lower bound for the growth of the fundamental group of a closed hyperbolic 3-manifold is a consequence of the following general result of Besson, Courtois and Gallot for Hadamard manifolds, see \cite{Besson1}, \cite{Besson2}.
\begin{theorem}[Besson-Courtois-Gallot]\label{Besso}
Let $(M^{n},g)$ be a Hadamard manifold whose sectional curvature satisfies $-a^{2}\leq K_{g}\leq -1$. Let $\Gamma$ be a discrete and finitely generated subgroup of ${\rm Isom}(M^{n}, g)$, then either $\Gamma$ is virtually nilpotent or its algebraic entropy is bounded below by a uniform explicit positive constant $c(n, a)$.
\end{theorem}
Recall that the algebraic entropy of a finitely generated group is defined as the logarithm of its uniform exponential growth
\begin{align}\notag
{\rm Ent}(\Gamma):= \log(\omega(G)).
\end{align}

Thus, let $M^{3}$ be a closed hyperbolic manifolds. It is well known that $\pi_{1}(M^{3})$ is a torsion free discrete cocompact subgroup of ${\rm Isom}(\setH^{3})$. Moreover it cannot be virtually nilpotent as follows for example from theorem \ref{Milnor}. Theorem \ref{Besso} then implies the following bound
\begin{align}\label{1}
\omega(\pi_{1}(M^{3}))\geq e^{c(3, 1)}> 1.
\end{align}
\begin{remark}
 Since the isometry group of $\setH^{3}$ is linear, a similar estimate could have been alternatively derived from the uniform Tits alternative for linear groups proved by Breuillard and Gelander in \cite{Breu}.
\end{remark}
\subsection*{Seifert fibred spaces}
We first describe which Seifert fibre spaces admit metrics of nonpositive curvature. There is a completely satisfactory division of closed Seifert fibre spaces into six classes according to which geometric structure they admit, see \cite{Scott} section 4. This geometric division allow us to conclude that the fundamental group of a closed Seifert fibre space modeled on $\setS^{3}$, $\setS^{2}\times \setR$, $\setE^{3}$, ${\rm Nil}^{3}$ must be finite, virtually infinite abelian or virtually nilpotent. Theorem \ref{Avez} then implies that these Seifert spaces cannot admit nonflat nonpositively curved metrics. For the sake of completeness we remark that the Bieberbach theorem implies
that the Seifert fibre spaces modeled on $\setS^{3}$, $\setS^{2}\times \setR$, ${\rm Nil}^{3}$ cannot admit flat metrics. It remains to study the Seifert fibre spaces which admit a geometric structure modeled on ${\rm \setH}^{2}\times \setR$ and ${\rm SL}(2,\setR)$. A Seifert fibre space modeled on $\setH^{2}\times \setR$ obviously admits nonpositively curved metric. Thus, let $M^{3}$ be a closed Seifert fibre bundle which admit a ${\rm SL}(2,\setR)$-geometry equipped with a nonpositively curved metric. Note that a ${\rm SL}(2,\setR)$-geometry cannot be nonpositively curved, see \cite{Milnor3}. Now, $ M^{3}$ is finitely covered by a orientable $\setS^{1}$-bundle over a closed orientable hyperbolic surface $F_{g}$. The fundamental group of such a bundle has a nontrivial $\setZ$-center, see \cite{Hempel}.
Therefore, by a theorem of Eberlein \cite{Eberlein}, it must be finitely covered by the trivial bundle $\setS^{1}\times F_{g}$. This implies that the Euler class of $M^{3}$ is zero which is a contradiction.

We can now study the exponential growth rate of the fundamental group of a closed Seifert fibre space $M^{3}$ modeled on ${\rm \setH}^{2}\times \setR$. Recall that $\pi_{1}(M^{3})$ fits in the short exact sequence
\begin{align}\notag
1\rightarrow\setZ\rightarrow\pi_{1}(M^{3})\rightarrow\pi_{1}(K)\rightarrow 1
\end{align}
where $K$ is the base orbifold. In this case $K$ is a closed good hyperbolic 2-orbifold, theorem \ref{Besso} then implies
\begin{align}\notag
\omega(\pi_{1}(K))\geq e^{c(2,1)}.
\end{align}
We conclude that
\begin{align}\label{2}
\omega(\pi_{1}(M^{3}))\geq\omega(\pi_{1}(K))\geq e^{c(2,1)}> 1.
\end{align}
\begin{remark}
Exactly the same uniform estimate holds for closed Seifert fibre spaces modeled on ${\rm SL}(2,\setR)$.
\end{remark}
\begin{remark}
If we restrict our attention to the class of $\setS^{1}$-bundle over surfaces of genus $g\geq 2$, the short exact sequence
\begin{align}\notag
1\rightarrow\setZ\rightarrow\pi_{1}(M^{3})\rightarrow\pi_{1}(F_{g})\rightarrow 1
\end{align}
provides the sharper estimate
\begin{align}\notag
\omega(\pi_{1}(M^{3}))\geq\omega(\pi_{1}(F_{g}))\geq 4g-3.
\end{align}
Nevertheless, since we cannot control the number of Seifert fibre spaces covered by these $\setS^{1}$-bundles there is no hope to obtain a uniform bound by simply using the Shalen-Wagreich lemma, see proposition 3.3 in \cite{Shalen}.
\end{remark}

\subsection*{Closed irreducible 3-manifolds with incompressible tori}
In what follows we will make use of the Bucher-de la Harpe estimates for the exponential growth of free products with amalgamation and HNN-extensions, for more details see \cite{Harpe} and \cite{Harpe1}.
\begin{theorem}[Bucher-de la Harpe]\label{Bucher}
Let $G=G_{1}\ast_{H}G_{2}$ be a finitely generated free product with amalgamation satisfying the condition $([G_{1}: H]-1)([G_{2}: H]-1)\geq 2$, then
\begin{align}\notag
\omega(G_{1}\ast_{H}G_{2})\geq \sqrt[4]{2}.
\end{align}
Let $H$, $K$ be two isomorphic subgroups of a finitely generated group $G$, and let $\theta: H\rightarrow K$ be an isomorphism. Consider the HNN-extension $G\ast^{\theta}_{H}$ and assume moreover that $[G: H]+[G: K]\geq 3$. Then
\begin{align}\notag
\omega(G\ast^{\theta}_{H})\geq \sqrt[4]{2}.
\end{align}
\end{theorem}
Let $M^{3}$ be orientable closed and irreducible with an incompressible embedded tori $T^{2}\hookrightarrow M^{3}$. By cutting $M^{3}$ along the incompressible embedded tori $T^{2}$ we obtain a orientable compact irreducible 3-manifold with boundary $M^{3}\vert\,T^{2}$ possibly nonconnected. Now, in order to apply theorem \ref{Bucher} we have to understand which orientable compact irreducible 3-manifold with boundary has fundamental group with a finite index $\setZ^{2}$.
The fundamental group of such a manifold would in particular be solvable. \\  In \cite{Evans}, Evans and Moser classify which solvable groups can occur as the fundamental group of a 3-manifold. More precisely, see theorem 4.2. in \cite{Evans}, they show that, up to connected sum with homotopy spheres, the only orientable irreducible compact 3-manifold with boundary and a finite index $\setZ^{2}$ are the torus times the interval or the twisted I-bundle over the Klein bottle.
Thus, if $\pi_{0}(M^{3}\vert\,T^{2})\neq 0$ and at least one of the connected components is not a twisted I-bundle over the Klein bottle we can conclude $\pi_{1}(M^{3})$ is a free product with $\setZ^{2}$-amalgamation that satisfies the hypothesis of theorem \ref{Bucher}. We then conclude that $\omega(\pi_{1}(M^{3}))\geq \sqrt[4]{2}$. Similarly, if $\pi_{0}(M^{3}\vert\,T^{2})= 0$ and $M^{3}\vert\,T^{2}$ is not the torus times the interval we conclude that $\pi_{1}(M^{3})$ is a HNN-extension as in theorem \ref{Bucher} and again the same bound holds. The remaining cases are now excluded if we assume that $M^{3}$ is nonpositively curved. In fact, if $M^{3}\vert\,T^{2}$ is the torus times the interval or if $M^{3}$ is the union of two twisted I-bundle over the Klein bottle, then it must be finitely covered by a torus bundle over the circle. It follows $\pi_{1}(M^{3})$ is virtually solvable. We can then conclude by the Gromoll-Wolf structure theorem for closed nonpositively curved manifolds, see \cite{Gromoll}. In summary, for a closed orientable irreducible 3-manifold $M^{3}$ with an incompressible tori, for which the fundamental group is not virtually solvable, we have
\begin{align}\label{3}
\omega(\pi_{1}(M^{3}))\geq \sqrt[4]{2}.
\end{align}

Combining \eqref{1}, \eqref{2}, \eqref{3} we obtain an explicit estimate for the constant $C$ in theorem \ref{Di Cerbo}.
\begin{remark}
The argument here proposed applies equally well to all irreducible closed 3-manifolds with nonempty JSJ decomposition and nonvirtually solvable fundamental groups. On the other hand there are example of irreducible graph manifolds, with nonvirtually solvable fundamental group, which do not admit nonpositively curved metrics, see \cite{Leeb1}. Nevertheless, if at least one of the component in the JSJ decomposition is hyperbolic then the manifold admits a nonpositively curved metric, see again \cite{Leeb1}.
\end{remark}

\section{On the growth of closed three manifold groups}\label{123}

The purpose of this section is to study the exponential growth rate of the remaining closed 3-manifold groups, that is, nonprime closed 3-manifolds and prime 3-manifolds with virtually solvable fundamental groups.
\subsection*{Nonprime closed 3-manifold groups}
We first note the following group theoretic result, due to Bucher \cite{Harpe}.
\begin{proposition}[Bucher]\label{Bucher2}
Let $G=G_{1}*G_{2}$ be a free product of two finitely generated groups with $|G_{1}|\geq 2$ and $|G_{2}|\geq 3$. Then $G$ is of uniform exponential growth and $\omega(G)\geq \sqrt{2}$.
\end{proposition}
Thus, let $M^{3}$ be a closed, nonprime, orientable 3-manifold. By the Kneser decomposition theorem $M^{3}$ has a nontrivial finite decomposition as a connected sum
\begin{displaymath}
M^{3}=M_{1}\#M_{2}\#...\#M_{j}\#\biggl(k\#(\setS^{2} \times \setS^{1})\biggr)
\end{displaymath}
where each of the $M_{i}$ is irreducible, see \cite{Milnor2}, \cite{Hempel}. The Seifert-Van Kampen theorem implies that $\pi_{1}(M^{3})$ can be expressed as the free product
\begin{align}\notag
\pi_{1}(M^{3})=\pi_{1}(M_{1})*...*\pi_{1}(M_{j})\underbrace{*\setZ*...*\setZ}_{k-times}.
\end{align}
Proposition \ref{Bucher2} then implies:
\begin{proposition}
Let $M^{3}$ be a closed, orientable, nonprime 3-manifolds such that $\pi_{1}(M^{3})\neq \setZ_{2}*\setZ_{2}$. Then  $\omega(\pi_{1}(M^{3}))\geq \sqrt{2}$.
\end{proposition}
In other words, any nonprime 3-manifold group of exponential growth is of uniform exponential growth bounded from below by a universal constant.

\subsection*{Virtually solvable closed 3-manifold groups}
It remains to study the growth of virtually solvable 3-manifold groups. Note that most of them are prime since the only nontrivial free product that is solvable is $\setZ_{2}*\setZ_{2}$. Moreover if we restrict our attention to virtually solvable 3-manifold groups of exponential growth, we then know that they must be finite extensions of the fundamental groups of some hyperbolic torus bundle over $\setS ^{1}$, see \cite{Evans} and \cite{Scott}. More precisely, we know that any 3-manifold with virtually solvable fundamental group of exponential growth admits a covering of degree at most four which is a hyperbolic torus bundle over $\setS^{1}$, see \cite{Scott}. Recall that a hyperbolic torus bundle over the circle is obtained from $T^{2}\times [0, 1]$ by identification of $(x, 1)$ with $(Ax, 0)$, where $A\in GL(2, \setZ)$ is a hyperbolic, i.e., none of its eigenvalues has absolute value one and $|det(A)|=1$. Now, the fundamental group of a hyperbolic torus bundle $M^{3}_{A}$ is a split extension of type $\setZ^{2}\rtimes_{A}\setZ$, as follows from the short exact sequence
\begin{align}\notag
1\rightarrow\setZ^{2}\rightarrow\pi_{1}(M^{3}_{A})\rightarrow\setZ\rightarrow 1.
\end{align}
The uniform exponential growth of this particular class of polycyclic groups is well understood, see \cite{Osin} and the bibliography there. In fact, given a group $G$ expressed as a split extension $\setZ^{n}\rtimes_{A}\setZ$ with $A$ a hyperbolic map, denoted with $Spec(A)$ the set of all eigenvalues of $A$, Osin \cite{Osin} shows:
\begin{align}\label{Osin1}
\omega(G)\geq 2^{\frac{\log\Lambda(A)}{\log(2)+\log\Lambda(A)}},
\end{align}
where
\begin{align}\notag
\Lambda(A):=\max_{\lambda \in Spec(A)}|\lambda|.
\end{align}
Now, an easy computation shows that in the case of hyperbolic torus bundle $M^{3}_{A}$ we have $\Lambda(A)> 2$. The estimate in \eqref{Osin1} then implies the bound $\omega(\pi_{1}(M^{3}_{A}))> 2^{\frac{1}{6}}$. Finally, since we can control uniformly the number of 3-manifolds covered by hyperbolic torus bundles, the Shalen-Wagreich lemma can be used to obtain a uniform lower bound for the exponential growth rate of virtually solvable 3-manifold groups. The following proposition is then established.
\begin{proposition}
There exists a universal constant $K>1$ such that for any virtually solvable closed 3-manifold group $G$ of exponential growth one has $\omega(G)\geq K$.
\end{proposition}

\section{Final remarks}

Combining the results above we obtain that a infinite closed 3-manifold group is either of polynomial growth or its uniform exponential growth rate is bounded from below by a universal constant. The maximal order of polynomial growth can also be derived. In fact, the explicit classification of nilpotent groups that arise as fundamental group of 3-manifolds, implies that the discrete 3-dimensional Heisenberg group provides the maximal order of growth. It is well known that this nilpotent group has polynomial growth of order four. This can be derived directly see \cite{Milnor1}, or by using the general Bass formula for the growth of finitely generated discrete nilpotent groups \cite{Harpe}.

We finally remark that analogous results hold in the case of closed good 3-orbifolds. The details will be presented elsewhere.

\section{Acknowledgements}
The author would like to thank G. Besson, G. Gallot, J. Lott for several useful discussions. He is also grateful to S. Maillot for invaluable advices and to his advisor C. LeBrun for constant support.

\end{document}